\documentclass[oneside,english]{amsart}
\usepackage[T1]{fontenc}
\usepackage[latin9]{inputenc}
\usepackage{amstext}
\usepackage{amsthm}
\usepackage{amssymb}

\makeatletter
\numberwithin{equation}{section}
\numberwithin{figure}{section}
\theoremstyle{plain}
\newtheorem{thm}{\protect\theoremname}
\theoremstyle{plain}
\newtheorem{lem}[thm]{\protect\lemmaname}
\newtheorem{prop}[thm]{Proposition}
\newtheorem{claim}[thm]{Claim}
\newtheorem{rem}[thm]{Remark}

\makeatother

\usepackage{babel}
\providecommand{\lemmaname}{Lemma}
\providecommand{\theoremname}{Theorem}
\providecommand{\theoremname}{Claim}
\providecommand{\theoremname}{Preposition}

\title{On $t$-Intersecting Families of Permutations}

\author{Nathan Keller}
\thanks{Department of Mathematics, Bar-Ilan University. \texttt{Nathan.Keller@biu.ac.il}. Supported by the Israel Science Foundation (grant no.~2669/21).}
\author{Noam Lifshitz}
\thanks{Einstein institute of Mathematics, Hebrew University. \texttt{noamlifshitz@gmail.com}. Supported by the Israel Science Foundation (grant no.~1980/22).}  
\author{Dor Minzer}
\thanks{Department of Mathematics, MIT. \texttt{minzer.dor@gmail.com}. Supported by a Sloan Research
Fellowship, NSF CCF award 2227876 and 
NSF CAREER award 2239160.} 
\author{Ohad Sheinfeld}
\thanks{Department of Mathematics, Bar-Ilan University. \texttt{oshenfeld@gmail.com}}

\begin{document}

\maketitle

\begin{abstract}
    We prove that there exists a constant $c_0$ such that for any $t \in \mathbb{N}$ and any $n\geq c_0 t$, if $A \subset S_n$ is a $t$-intersecting family of permutations then$|A|\leq (n-t)!$. Furthermore, if $|A|\ge 0.75(n-t)!$ then there exist $i_1,\ldots,i_t$ and $j_1,\ldots,j_t$ such that $\sigma(i_1)=j_1,\ldots,\sigma(i_t)=j_t$ holds for any $\sigma \in A$. This shows that the conjectures of Deza and Frankl (1977) and of Cameron (1988) on $t$-intersecting families of permutations hold for all $t \leq c_0 n$. Our proof method, based on hypercontractivity for global functions, does not use the specific structure of permutations, and applies in general to $t$-intersecting sub-families of `pseudorandom' families in $\{1,2,\ldots,n\}^n$, like $S_n$.
\end{abstract}

\section{Introduction}

\subsection{Background} A family $F$ of subsets of $[n]=\{1,2,\ldots,n\}$ is called \emph{$t$-intersecting} if for any $A,B \in F$, we have $|A \cap B| \geq t$. For $t=1$, such families are simply called `intersecting'. In 1961, Erd\H{o}s, Ko, and Rado~\cite{EKR61} 
proved that for any $k< n/2$, the maximal size of an intersecting family of $k$-subsets of $[n]$ is ${n-1}\choose{k-1}$, and the maximum is obtained only for \emph{dictatorships}, i.e., families of the form $\{S \subset [n]: |S|=k, i \in S\}$, for some $i \in [n]$. The Erd\H{o}s-Ko-Rado theorem launched a field of research in extremal combinatorics which studies families of finite sets with various restrictions on intersections between the sets in the family (see the survey~\cite{FT16}). One of the central problems in this field is determining the maximal size of a $t$-intersecting family $F \subset \mathcal{U}$, for various `universes' $\mathcal{U}$. 

In the basic setting, where $\mathcal{U}$ consists of all $k$-element subsets of $[n]$, Erd\H{o}s, Ko, and Rado showed in~\cite{EKR61} that for $n \geq n_0(k,t)$, the maximal size is ${n-t}\choose{k-t}$, and asked, what is the minimal number $n_0(k,t)$ for which this upper bound holds. This question was solved by Frankl~\cite{F78} for all $t \geq 15$, and then by Wilson~\cite{W84} for all $t$: they showed that the minimal number is $n_0(k,t)=(k-t+1)(t+1)$. Furthermore, for all $n>n_0(k,t)$, the maximal size is obtained only for \emph{$t$-umvirates}, i.e., families of the form $\{S \subset [n]: |S|=k, T \subset S\}$, for some $|T|=t$. Friedgut~\cite{Fri08} and Ellis, Keller, and Lifshitz~\cite{EKL19} obtained \emph{stability versions} of the theorem, which assert that if the size of a $t$-intesecting family is `close' to the maximum possible size, then it is essentially contained in a $t$-umvirate. For $t=1$ (i.e., intersecting families), stability versions were obtained much earlier by Hilton and Milner~\cite{HM67} and by Frankl~\cite{Fra87}. The more general question of determining the maximal possible size of a $t$-intersecting family for any triple $(n,k,t)$ was solved in a celebrated theorem of Ahlswede and Khachatrian~\cite{AK97}, which played a central role in applications of intersection theorems to computer science (see, e.g.,~\cite{DS05}).

Extremal problems on $t$-intersecting families were studied in various other `universes' $\mathcal{U}$ -- e.g., for graphs~\cite{EFF12}, set partitions~\cite{MM05}, and linear maps~\cite{Linear-maps}. Arguably, the most thoroughly studied `universe' (except for $k$-subsets of $[n]$) is $t$-intersecting families of permutations, i.e., families $F \subset S_n$ such that for any $\sigma,\tau \in F$, there exist $i_1,\ldots,i_t$ with $\sigma(i_j)=\tau(i_j)$ for $j=1,\ldots,t$.

In 1977, Deza and Frankl~\cite{DF77} proved an analogue of the Erd\H{o}s-Ko-Rado theorem for permutations: they showed that the maximal size of an intersecting family $F \subset S_n$ is $(n-1)!$, which is obtained for the \emph{dictatorship} families $(S_n)_{i \to j} = \{\sigma \in S_n: \sigma(i)=j\}$. Showing that dictatorships are the only maximal families, turned out to be harder. This was conjectured by Deza and Frankl~\cite{DF77}, but was showed only in 2003 by Cameron and Ku~\cite{CK03} and (independently) by Larose and Malvenuto~\cite{LM04}. For a general $t$, Deza and Frankl conjectured (in the same paper~\cite{DF77}) that for all $n \geq n_0(t)$, the maximal size is $(n-t)!$, which is obtained for the \emph{$t$-umvirate} families 
\[
(S_n)_{i_1 \to j_1,\ldots,i_t \to j_t} = \{\sigma \in S_n: \forall 1 \leq j \leq t, \sigma(i_1)=j_t\}. 
\]
Cameron~\cite{Cam88} conjectured in 1986 that the $t$-umvirates are the only maximum-sized families for all $n \geq n_0(t)$. 

The conjectures were open for several decades. Only in 2011, Ellis, Friedgut, and Pilpel~\cite{EFP11} proved that the Deza-Frankl conjecture holds when $n$ is at least double exponentially large in $t$ (equivalently, when $t=O(\log \log n)$). Ellis et al.~\cite{EFP11} also showed a decoupled version which asserts that if $F_1,F_2 \subset S_n$ are \emph{cross $t$-intersecting}, meaning that for any $\sigma \in F_1$ and $\tau \in F_2$ there exist $i_1,\ldots,i_t$ with $\sigma(i_j)=\tau(i_j)$ for $j=1,\ldots,t$, then $|F_1||F_2|\leq (n-t)!^2$, proving (in a strong form) a conjecture of Leader~\cite{Lea05}. The proof of Ellis et al.~is quite involved and uses Fourier analysis on the symmetric group which entails representation theory. Essentially, their idea was to design a Cayley graph over $S_n$ whose adjacency matrix has a specific set of eigenvalues and all of whose independent sets are $t$-intersecting sets of permutations. They then employed the Hoffman bound, which upper bounds independent sets in graphs in terms of their minimal eigenvalue to obtain their result. Roughly at the same time, Ellis~\cite{Ell11} proved that Cameron's conjecture holds when $n$ is at least double exponentially large in $t$. He also showed a stability version, which asserts that if $|F|>(1-1/e+o(1))(n-t)!$, then $F$ must be contained in a $t$-umvirate. 

As a natural `next step', Ellis et al.~\cite{EFP11} conjectured that an analogue of the Ahlswede-Khachatrian theorem holds for $t$-intersecting families of permutations. In particular, they conjectured that for all $t<n/2$, the $t$-umvirates are the unique maximum-sized $t$-intersecting families in $S_n$ -- namely, that the aforementioned conjectures of Deza and Frankl~\cite{DF77} and of Cameron~\cite{Cam88} hold for all $n>2t$. Note that this bound on $t$ is the maximal possible, as for $t=n/2$, the family $F\subset S_n$ of all permutations that have at least $t+1$ fixed points among $\{1,2,\ldots,t+2\}$ is $t$-intersecting and its size is $(t+2)\cdot(n-t-1)!-(t+1)\cdot(n-t-2)!$, which is larger that $(n-t)!$ for all $n\ge 8$.

The first step toward proving the conjecture of Ellis et al.~was obtained by Ellis and Lifshitz~\cite{DN22} who showed that the conjecture holds for all $t=O(\frac{\log n}{\log\log n})$. In fact, they showed a stronger result: If no $\sigma,\tau \in F$ agree on exactly $t-1$ elements, then $|F|\leq (n-t)!$, and the $t$-umvirates are the only maximum-sized families. This setting, called `forbidding one intersection', is considerably harder than the $t$-intersecting setting, and is not fully resolved even for families of $k$-subsets of $[n]$ (see~\cite{EKL16} for the most recent results on this well-known problem, called the Erd\H{o}s-S\'{o}s problem). The proof of Ellis and Lifshitz uses the discrete Fourier-analytic \emph{junta method}~\cite{NN21}, along with a representation-theoretic argument. The main idea was to employ a dichotomy between a $t$-umvirate-structure and a pseudorandomness notion called quasiregularity/globalness. They partitioned the $t$-intersecting families into `global' families. They then proved structural results about the Fourier decomposition of the indicators of global sets, and used these to obtain their bound. 

A major breakthrough was obtained in a recent work of Kupavskii and Zakharov~\cite{KZ22} who proved that the conjectures of Deza and Frankl~\cite{DF77} and of Cameron \cite{Cam88} hold for all $t=O(\frac{n}{(\log(n))^2})$, along with a stability version. Kupavskii and Zakharov proved the same statement also in the harder `forbidding one intersection' setting, for all $t=\tilde{O}(n^{1/3})$. Unlike all previous works in this direction (e.g.,~\cite{CK03,EFP11,DN22,DF77,LM04}), the work of Kupavskii and Zakharov does not use any specific properties of $S_n$ (and in particular, its representation theory). Instead, they prove that the  conjecture holds for the general setting of $t$-intersecting families in $\mathcal{U}$, where the `universe' $\mathcal{U}$ is a `pseudorandom' sub-family of ${N}\choose{M}$, and show that $S_n$ can be viewed as a `pseudorandom' sub-family of ${n^2}\choose{n}$. The technique of Kupavskii and Zakharov is called the \emph{spread approximations method}. Similarly to \cite{DN22}, Kupavskii and Zakharov reduce their analysis to the study of global sets. Their method relies on the recent breakthrough results of Alweiss, Lovett, Wu, and Zhang~\cite{ALWZ21} on the sunflower problem in place of the Fourier analytic approach of \cite{DN22}. 

Much less is known when $t$ may be as large as linear in $n$. The only result in this setting so far was obtained by Keevash and Long~\cite{KL17} who proved that for any $\epsilon > 0$, there exists $\delta>0$ such that if $\epsilon \cdot n < t < (1-\epsilon)\cdot n$ and $F \subset S_n$ does not contain two permutations that agree on exactly $t-1$ points, then $|F| < (n!)^{1-\delta}$.

\subsection{Our results} We prove that the conjectures of Deza and Frankl~\cite{DF77} and of Cameron~\cite{Cam88} hold for all $t \leq c_0 n$, for a universal constant $c_0$. Furthermore, we prove a decoupled version and a stability version. Our main result is the following:
\begin{thm}\label{Thm:Main}
    There exists $c_0>0$ such that the following holds for all $t \in \mathbb{N}$ and all $n\geq c_0 t$. Let $A,B \subset S_n$ be cross $t$-intersecting families. Then: 
    \begin{itemize}
        \item $|A||B| \leq (n-t)!^2$.

        \item If $|A||B| \geq 0.75 (n-t)!^2$, then there exist $i_1,\ldots,i_t$ and $j_1,\ldots,j_t$ such that $A,B \subset (S_n)_{i_1\to j_1,\ldots,i_t \to j_t}$.

        \item If $A=B$ and $|A| \geq 0.75 (n-t)!$, then there exist $i_1,\ldots,i_t$ and $j_1,\ldots,j_t$ such that $A\subset (S_n)_{i_1\to j_1,\ldots,i_t \to j_t}$.
    \end{itemize}
\end{thm}
The constant $0.75$ was chosen for convenience of the proof, and can be replaced with any constant larger than $1-1/e$ (see Remark~\ref{Rem:Optimal-constant}). As was shown by Ellis~\cite{Ell11} for $t=1$, the constant $1-1/e$ is optimal. Indeed, if $\sigma \in S_n$ is the permutation which interchanges $1$ and $n$ and leaves all other points fixed, $A=(S_n)_{1\to 1,\ldots,t\to t} \cup \{\sigma\}$, and $B=\{\tau \in (S_n)_{1\to 1,\ldots,t\to t} : |\tau \cap \sigma|\geq t\}$, then $A,B$ are cross $t$-intersecting and $|A||B|=(1-1/e-o(1))(n-t)!^2$. Similarly, $B \cup \{\sigma\}$ is $t$-intersecting and $|B|=(1-1/e-o(1))(n-t)!$, which shows that the constant $0.75$ in the second and the third assertions of the theorem cannot be replaced by any constant smaller than $1-1/e$.

Our methods can also be used to prove a variant of Theorem~\ref{Thm:Main} in the harder `forbidding one intersection' setting, albeit with a smaller upper bound on $t$:
\begin{prop}
    Let $t,n \in \mathbb{N}$, such that $t\leq c'_0\sqrt{n}/\log{n}$, where $c'_0$ is an absolute constant. Let $A,B \subset S_n$ such that for any $\sigma \in A, \tau \in B$, $|\sigma \cap \tau|\neq t-1$. Then $|A||B| \leq (n-t)!^2$.
\end{prop}
As this result is more technical, it will be presented in a future paper~\cite{KLS23+}.

\subsection{Our techniques and outline of the proof} Like in the work of Kupavskii and Zakharov~\cite{KZ22} we show that the assertion of the theorem holds for cross $t$-intersecting families in any `pseudorandom' subfamily of a host space. In their case the host space is $\binom{[n^2]}{n}$ and in our case the host space is  $(\mathbb{Z}/n\mathbb{Z})^n$. The method is in the spirit of the Green and Tao~\cite{GT08} philosophy that large subsets of pseudorandom subsets of a space $X$ behave like large subsets of $X$. Green and Tao viewed the primes as a large subset of the pseudorandom `almost primes'. Here we view $S_n$ as a pseudorandom subset of the Abelian group $(\mathbb{Z}/n\mathbb{Z})^n.$ This allows us to employ an analytical tool from the theory of product spaces known as hypercontractivity for global functions developed by Keevash et al. \cite{KLLM21}.

Our proof consists of two main steps:
\begin{enumerate}
    \item \emph{Cross $t$-intersecting families have a density bump inside a dictatorship.} We show that if $A,B \subset S_n$ are cross $t$-intersecting then there exist $i,j$ such that the relative density of $A$ and $B$ inside the dictatorship $(S_n)_{i\to j}$ is significantly larger than the density of $A$ and $B$, respectively.

    \item \emph{Upgrading density bump inside a dictatorship into containment in a $t$-umvirate.} We show that the first step can be applied sequentially to deduce that $A,B$ are essentially contained in a $t$-umvirate $(S_n)_{i_1\to j_1,\ldots,i_t \to j_t}$ for some $i_1,\ldots,i_t$ and $j_1,\ldots,j_t$. 
\end{enumerate}
The second step, presented in Section~\ref{sec:main}, is an inductive argument which requires incorporating a stability statement as part of the proof.

The first step, which encapsulates the `pseudorandomness vs. structure' philosophy, is the more involved one. Here, we first use a coupling technique 
to embed our families into $\{0,1\}^{n^2}$, endowed with the biased measure $\mu_p$ (for a value of $p$ that depends on $t$), and argue in $(\{0,1\}^{n^2},\mu_p)$. On the contra-positive, we assume that $A$ and $B$ (viewed as subfamilies of $\{0,1\}^{n^2}$) don't have a density bump inside a dictatorship (which is our notion of pseudorandomness), and reach a contradiction. 

Our first sub-step is to upgrade this weak pseudorandomness notion to the stronger notion of \emph{globalness} presented in~\cite{KLLM21}. This step, presented in Section~\ref{sec:density-bump}, proceeds by constructing restrictions of $A$ and $B$ which are sufficiently large, global, and cross $t$-intersecting. 

The second sub-step, presented in Section~\ref{sec:global-t-intersecting-are-small}, allows using the restrictions of $A$ and $B$ to obtain a contradiction. We observe that the $\mu_{1/3}$ measures of $A$ and $B$ cannot be `too large', since by the FKG correlation inequality~\cite{FKG}, this would imply that $A\cap B$ is a `too large' $t$-intersecting family w.r.t.~to $\mu_{1/3}$, contradicting the biased-measure version of the Ahlswede-Khachatrian theorem~\cite{Filmus17}. We then complete the proof by presenting a general \emph{sharp threshold} argument which shows that if the $\mu_p$ measure of a global subfamily of $\{0,1\}^{n^2}$ is not extremely small, then its $\mu_{\frac{1}{3}}$ measure is large. This sharp threshold argument relies on the \emph{level-$d$ inequality for global functions}, recently proved in~\cite{Omri}.

\section{Preliminaries}
\label{sec:Preliminaries}


\subsection{Notation}
\emph{Constants.} We use many different constants throughout the paper. The constants whose values do matter are numbered (e.g., $c_0,\ldots,c_3,n_0$), while the values of all other constants (like $C,C',g,G$) do not matter, and such notations are re-used between different proofs. All constants are positive.

\emph{The `measure' of functions.} For a measure space $(K,\mu)$ and a function $f:K\to \{0,1\}$, we denote $\mu(f):=\Pr_{\mu}[f(x)=1]$.



\subsection{The biased measure on the discrete cube}
For $0<p<1$, the biased measure on $\{0,1\}^n$ is defined by $\mu_p(x)=p^{\sum x_i} (1-p)^{n-\sum x_i}$. We shall use several classical notions and results regarding functions in $L^2(\{0,1\}^n,\mu_p)$.

\subsubsection*{Fourier expansion.} For each $i \in [n]$, let $\chi_i(x)=\sqrt{\frac{1-p}{p}}$ if $x_i=1$ and $\chi_i(x)=-\sqrt{\frac{p}{1-p}}$ if $x_i=0$. For each $S \subset [n]$, let $\chi_S=\prod_{i \in S}\chi_i$. The set of functions $\{\chi_S\}_{S\subset [n]}$ is an orthonormal basis of $L^2(\{0,1\}^n,\mu_p)$. Hence, any function $f \in L^2(\{0,1\}^n,\mu_p)$ admits a unique representation $f=\sum_S \hat f(S)\chi_S$. The coefficients $\hat f(S)$ are called the Fourier coefficients of $f$. The \emph{level} of the coefficient $\hat f(S)$ is $|S|$, and we denote by $f^{=d}$ the $d$'th level of $f$, namely, $f^{=d}=\sum_{|S|=d} \hat f(S)\chi_S$.

In cases where the underlying measure $\mu_p$ is not clear from the context, we denote the characters by $\chi^p_S$.

\subsubsection*{$(q,p)$-biased distribution.} For $0<q<p<1$, the $(q,p)$-biased distribution $D(q,p)$ is the unique probability distribution on pairs $(x,y) \in \{0,1\}^n \times \{0,1\}^n$ such that:
\begin{enumerate}
    \item All the pairs $(x_i,y_i)$ are mutually independent.
    \item For all $i$, we have $\Pr[x_i=1]=q$ and $\Pr[y_i=1]=p$.
    \item For all $i$, we have $x_i \leq y_i$ with probability $1$.
\end{enumerate}

\subsubsection*{One-sided noise operator.}
For $0<q<p<1$, the one-sided noise operator $T_{q \rightarrow p}:L^2(\{0,1\}^n,\mu_q) \rightarrow L^2(\{0,1\}^n,\mu_p)$ is defined by $T_{q \rightarrow p}f(y)=\mathbb{E}[f(x)]$, where $(x,y)$ is distributed according to $D(q,p)$. 

This operator, first introduced in~\cite{ABGM14}, is called a `one-sided noise operator', since $x$ is a `noisy' version of $y$ obtained by (possibly) changing some coordinates from 1 to 0, and leaving the zero coordinates of $y$ unchanged. (This lies in contrast with the standard noise operator, in which coordinates may be changed from 0 to 1 as well). 

As was shown in~\cite[Lemma~1]{Lifshitz20}, the one-sided noise operator admits the following convenient Fourier representation. For $f=\sum_S \hat f(S)\chi^q_S$,
\begin{equation}\label{Eq:One-sided noise}
    T_{q \rightarrow p}f = \sum_S \rho^{|S|} \hat f(S) \chi^p_S,
\end{equation}
where $\rho=\sqrt{\frac{q(1-p)}{p(1-q)}}$.

\subsubsection*{The FKG correlation inequality.}
We naturally identify elements of $\{0,1\}^n$ with subsets of $[n]$. The $p$-biased measure of a family $F \subset \{0,1\}^n$ is $\sum_{A \in F} \mu_p(A)$. A family $F$ is called \emph{monotone} if $(A \in F) \wedge (A \subset B) \Rightarrow (B \in F)$. 

A special case of the classical FKG inequality~\cite{FKG} asserts that any two monotone families are non-negatively correlated with respect to any product measure on $\{0,1\}^n$. In particular, if $F,G \subset \{0,1\}^n$ are monotone families then for any $0<p<1$,
\begin{equation}\label{Eq:FKG}
    \mu_p(F \cap G) \geq \mu_p(F) \mu_p(G).
\end{equation}

\subsubsection*{The Ahlswede-Khachatrian theorem.}
We shall use the biased-measure version of the classical Ahlswede-Khachatrian theorem~\cite{AK97}, proved by Filmus in~\cite[Theorem~3.1,Corollary~3.2]{Filmus17}.
\begin{thm}\label{Thm:Biased-AK}
 Let $F \subset \{0,1\}^n$ be a $t$-intersecting family. Then for any $r \geq 0$, if $\frac{r}{t+2r-1}<p<\frac{r+1}{t+2r+1}$ then $\mu_p(F) \leq \mu_p(\mathcal{F}_{t,r})$, where $\mathcal{F}_{t,r}=\{S:|S \cap [t+2r]| \geq t+r\}$.
\end{thm}
In particular, for $p=1/3$ (which is the case we will use), the theorem implies the upper bound $\mu_p(\mathcal{F}) \leq 0.85^t$.

\subsection{Restrictions and globalness}
Let $(K^n,\mu^n)$ be a product measure space, and let $f:K^n \rightarrow \{0,1\}$. For $S \subset [n]$ and $x \in K^{S}$, the \emph{restriction} of $f$ obtained by fixing the coordinates in $S$ to $x$ is denoted by $f_{S \rightarrow x}:K^{[n]\setminus S} \rightarrow \{0,1\}$. 

For $g>0$, a function $f$ is called \emph{$g$-global} if for any $S \subset [n]$ and any $x \in K^{S}$, we have $\mu^{n-|S|}(f_{S \rightarrow x}) \leq g^{|S|} \mu^{n}(f)$. This means that no fixing of coordinates increases the relative density of $f$ significantly.

\subsubsection*{$g$-global restrictions.} 
For any function $f: K^n \rightarrow \{0,1\}$ and any $g>0$, we may construct a $g$-global restriction of $f$ by choosing $S,x$ such that $\frac{\mu^{n-|S|}(f_{S\rightarrow x})}{g^{|S|}}$ is maximal over all the choices of $S,x$. (If there are several maxima, we pick one of them arbitrarily.) It is easy to see that $f_{S\rightarrow x}$ is indeed $g$-global, and also satisfies 
\begin{equation}\label{Eq:g-global}
    \mu^{n-|S|}(f_{S\rightarrow x}) \geq g^{|S|}\cdot \mu^n(f),
\end{equation} 
a fact which we will use several times. When $K=\mu_p$ and $\frac{1}{1-p} < g$, it's easy to see that the restriction $S\rightarrow x$ must be $S\rightarrow 1$, meaning the restriction is obtained by setting some coordinates to 1.

\subsubsection*{$d$-level inequality for global functions.} We will use the following sharp version of the $d$-level inequality for global functions on $(\{0,1\}^n,\mu_p)$, obtained in~\cite[Thm.~1.5]{Omri}. 
\begin{thm}\label{Thm:level-d}
    There exist $c_1,c_2>0$ such that the following holds. For any $g>0$, $n \in \mathbb{N}$, $0<p<1$, and any $g$-global function $f \colon \{0,1\}^n \to \{0,1\}$, we have:
    \[
    \|f^{= d}\|_2^2 \leq \mu_p(f)^2 \frac{c_2^d g^{2d} \log^d(\frac{1}{\mu_p(f)})}{d^d},
    \]
    for all $d \leq c_1\log(1/\mu_p(f))$. 
\end{thm}

\section{Global $t$-intersecting Families are Small}
\label{sec:global-t-intersecting-are-small}

The extremal families in the classical theorems on $t$-intersecting families, like the families $\mathcal{F}_{t,r}$ in the biased Ahlswede-Khachatrian theorem (Theorem~\ref{Thm:Biased-AK} above) are very asymmetric: membership in these families is determined by only a constant number of coordinates. Thus, it is natural to anticipate that imposing a `regularity' condition on the family will lead to a much stronger upper bound on its possible size. Results in this spirit were indeed proved in several settings. Ellis, Kalai and Narayanan~\cite{EKN20} showed that whenever $k \leq n/2- \omega(n)\log n/n$ for any function $\omega(n)$ that grows to infinity as $n \to \infty$, the maximum size of a symmetric intersecting family of $k$-element subsets of $[n]$ is $o({{n-1}\choose{k-1}})$ (while the maximum size of a general intersecting family of $k$-element subsets is exactly ${{n-1}\choose{k-1}}$ by the Erd\H{o}s-Ko-Rado theorem). They also obtained a similar result for symmetric families of subsets of $\{0,1\}^n$ w.r.t. the biased measure $\mu_p$, for $p \leq 1/2-\omega(n)\cdot \log n/n$. Ihringer and Kupavskii~\cite{IK19} obtained results in this spirit for regular families of $k$-element subsets of $[n]$ (i.e., families in which each element is included in the same number of sets); those results were recently quantitatively improved by Kupavskii and Zakharov~\cite{KZ22}. Eberhard, Kahn, Narayanan, and Spirkl~\cite{EKNS21} obtained similar results for symmetric intersecting families of vectors in $[m]^n$. 

In this section, we show that for any $t$ and any $p \leq p_0$ for a constant $p_0$, adding a globality assumption (which is much weaker than a regularity assumption) is already sufficient for deducing that a $t$-intersecting family has a small measure in $(\{0,1\}^n,\mu_p)$.

As this will be convenient for us in the sequel, we prove a decoupled version of the statement.
\begin{prop}\label{Prop:Global-cross-intersecting-are-small}
There exist constants $c_3,p_0>0$ such that the following holds for any $g>0$, $p \leq p_0$, $n \in \mathbb{N}$ and $t \in [n]$. Let $A,B\subseteq \{0,1\}^n$ be $g$-global, monotone, and cross $t$-intersecting. Then 
\[
\min(\mu_p(A),\mu_p(B)) \leq e^{-c_3 \cdot t \cdot \frac{1}{p} \cdot \frac{1}{g^2}}
\]
\end{prop}
To prove the proposition, we need the following sharp-threshold lemma for global functions:
\begin{lem}\label{Lemma:Sharp-threshold-global}
There exist constants $c_3,p_0>0$ such that the following holds for any $g>0$, $p \leq p_0$, $n \in \mathbb{N}$ and $t \in [n]$. Let $A\subseteq \{0,1\}^n$ be $g$-global and monotone. If $\mu_p(A)>e^{-c_3 \cdot t \cdot \frac{1}{p} \cdot \frac{1}{g^2}}$, then $\mu_{1/3}(A) \ge 0.99^t$.
\end{lem}

%

\begin{proof}
Throughout the proof we omit floors and ceilings for brevity. 
Let $f = 1_A \in L^2(\{0,1\}^n,\mu_p)$ and $g=1_A \in L^2(\{0,1\}^n,\mu_{1/3})$, and let $T_{p \to 1/3}$ be the one-sided noise operator defined in Section~\ref{sec:Preliminaries}. Recall that by~\eqref{Eq:One-sided noise}, if the Fourier expansion of $f$ w.r.t. $\mu_p$ is $f = \sum_S \hat f(S) \chi^p_S$ then
\[
T_{p \rightarrow 1/3}f = \sum_S \rho^{|S|} \hat f(S) \chi^{1/3}_S,
\]
where $\rho=\sqrt{\frac{(1-1/3)p}{1/3(1-p)}} \leq 2\sqrt{p}$.


Observe that 
\begin{equation}\label{Eq:Inner-product1}
\langle T_{p \rightarrow \frac{1}{3}}f, 1-g \rangle  = \mathbb{E}_{(x,y) \in D(p,1/3)}[f(x)(1-g(y))]
\end{equation}
where $D(p,1/3)$ is the $(p,1/3)$-biased distribution in which $x \sim \mu_p$, $y \sim \mu_{1/3}$, and $x \leq y$.
As $A$ is a monotone family, we have $f(x)=0$ whenever $g(y) = 0$, and thus, the right hand side of~\eqref{Eq:Inner-product1} is zero. 

Now, assume on the contrary that $\mu_\frac{1}{3}(g) <0.99^t$. We have:
\begin{align*}
\begin{split}
0 &= \langle T_{p \rightarrow \frac{1}{3}}f, 1-g \rangle = \sum_S \rho ^{|S|} \hat{f}(S) \widehat{(1-g)}(S) \\
&= \mu_p(f)\mu_\frac{1}{3}(1-g) - \sum_{|S|\ge 1} \rho ^{|S|} \hat{f}(S) \hat{g}(S),    
\end{split}
\end{align*}
and hence,
\begin{equation}\label{Eq:Product-of-measures1}
\mu_p(f)\cdot \mu_\frac{1}{3}(1-g) = \sum_{|S|\ge 1} \rho ^{|S|} \hat{f}(S) \hat{g}(S).
\end{equation}
We bound from above the absolute value of the RHS. By the Cauchy-Schwarz inequality,
\begin{align}\label{Eq:Product-of-measures2}
\begin{split}
    |\sum_{|S|\ge 1} &\rho ^{|S|} \hat{f}(S) \hat{g}(S)| \le \sum_{d=1}^{\infty} \rho^d \|f^{= d}\|_2 \cdot \|g^{= d}\|_2 \\ 
    &= \sum_{d=1}^{c_1\log(1/\mu_p(f))} \rho^d \|f^{= d}\|_2 \cdot \|g^{= d}\|_2 + \sum_{d=c_1\log(1/\mu_p(f))}^{\infty} \rho^d \|f^{= d}\|_2 \cdot \|g^{= d}\|_2, 
\end{split}
\end{align}
where $c_1$ is the constant from Theorem~\ref{Thm:level-d}.

Let us bound the first summand first. By assumption, we have $\mu_p(f)>e^{-c_3 \cdot t \cdot \frac{1}{p} \cdot \frac{1}{g^2}}$, or equivalently, $\log(\frac{1}{\mu_p(f)}) < \frac{c_3 \cdot t}{p\cdot g^2}$. Hence, by the level-$d$ inequality for global functions (i.e., Theorem~\ref{Thm:level-d} above), for all $d < c_1\log(1/\mu_p(f))$, we have 
\[
\|f^{= d}\|_2 \le \mu_p(f) \cdot \sqrt{\frac{c_2^d\cdot g^{2d} \cdot (\frac{c_3 \cdot t}{p\cdot g^2})^d}{d^d}} = \mu_p(f) \cdot \sqrt{\frac{c_2^d c_3^d \cdot t^d}{p^d d^d}},
\]
where $c_2$ is the constant from Theorem~\ref{Thm:level-d}.
Since $\rho \leq 2\sqrt{p}$, we have
\begin{align*}
\begin{split}
    \sum_{d=1}^{c_1\log(1/\mu_p(f))} \rho^d \|f^{= d}\|_2 \cdot \|g^{= d}\|_2 &\le \sum_{d=1}^{c_1\log(1/\mu_p(f))} \mu_p(f) \cdot\|g^{= d}\|_2 \cdot \left(\sqrt{\frac{4c_2 c_3 \cdot p\cdot t}{p\cdot d}}\right)^d \\ 
    &\le \sum_{d=1}^{c_1\log(1/\mu_p(f))} \mu_p(f) \cdot\|g^{= d}\|_2 \cdot \frac{1}{1000}\left(\sqrt{\frac{C c_3\cdot t}{d}}\right)^d,
\end{split}
\end{align*}
where $C$ (which may depend on $c_2$ but not on $c_3$) is a sufficiently large constant such that the last inequality holds. 

We further sub-divide the sum into three sub-sums:
\begin{itemize}
    \item $1 \leq d \leq C \cdot c_3 \cdot t$. In this range, we have $(\sqrt{\frac{C \cdot c_3 \cdot t}{d}})^d \le (\frac{C \cdot c_3\cdot t}{d})^d \le \frac{(C \cdot c_3 \cdot t)^d}{d!} $. Thus, we have
\begin{align*}
\begin{split}
    \sum_{d=1}^{C \cdot c_3\cdot t} \mu_p(f) \|g^{= d}\|_2 \cdot &\frac{1}{1000} \left(\sqrt{\frac{C \cdot c_3\cdot t}{d}}\right)^d \le \mu_p(f)  \|g\|_2 \sum_{d=1}^{C \cdot c_3\cdot t} \frac{1}{1000}\frac{(C \cdot c_3\cdot t)^d}{d!} \\
    &\le \mu_p(f) \cdot 0.995^t \cdot \frac{e^{C \cdot c_3 \cdot t}}{1000} \le 0.001 \mu_p(f), 
\end{split}
\end{align*}
where the penultimate inequality holds since $\|g\|_2 = (\mu_{1/3}(g))^{1/2}< 0.995^t$, and the last inequality holds, provided that $c_3$ is smaller than some constant that depends only on $C$. 

\item $C \cdot c_3\cdot t \le d \le 4C \cdot c_3\cdot t$. In this range, we have  $(\sqrt{\frac{C \cdot c_3\cdot t}{d}})^d \le 1$, and hence, 
\[
\sum_{C \cdot c_3\cdot t+1}^{4C \cdot c_3\cdot t} \mu_p(f) \|g^{= d}\|_2 \cdot \frac{1}{1000} (\sqrt{\frac{C \cdot c_3\cdot t}{d}})^d \le \mu_p(f) \cdot  0.995^t \cdot 3C \cdot c_3\cdot t \le 0.001 \mu_p(f),
\]
where the last inequality holds, provided that $c_3$ is smaller than some constant that depends only on $C$. 

\item $4C \cdot c_3\cdot t \le d \le c_1\log(1/\mu_p(f))$. In this range, we have $(\sqrt{\frac{C \cdot c_3\cdot t}{d}})^d \le (\frac{1}{2})^d$, and thus, 
\[
\sum_{d=4C \cdot c_3\cdot t}^{c_1\log(1/\mu_p(f))} \mu_p(f) \|g^{= d}\|_2 \cdot \frac{1}{1000}(\sqrt{\frac{C \cdot c_3\cdot t}{d}})^d \le \mu_p(f) \cdot \frac{1}{1000} \sum_{d=1}^{\infty} (\frac{1}{2})^d \le 0.001 \mu_p(f).
\]
\end{itemize}
Combining the sub-sums together, we obtain 
\begin{equation}\label{Eq:Product-of-measures3}
\sum_{d=1}^{c_1\log(1/\mu_p(f))} \rho^d \|f^{= d}\|_2 \cdot \|g^{= d}\|_2 \leq 0.003 \mu_p(f).
\end{equation}

Now, we bound from above the second summand in~\eqref{Eq:Product-of-measures2}. We have
\begin{align*}
    \sum_{d=c_1\log(1/\mu_p(f))}^{\infty} \rho^d \|f^{= d}\|_2 \cdot \|g^{= d}\|_2 \le  \sum_{d=c_1\log(1/\mu_p(f))}^{\infty} \rho^d \le 2\cdot \rho^{c_1\log(1/\mu_p(f))},
\end{align*}
where the last inequality holds provided that $p$ is sufficiently small. (As $\rho \leq 2\sqrt{p}$, it is sufficient that $p \leq 1/16$). 

We claim that 
\begin{equation}\label{Eq:Product-of-measures4}
    \sum_{d=c_1\log(1/\mu_p(f))}^{\infty} \rho^d \|f^{= d}\|_2 \cdot \|g^{= d}\|_2 \le 2\cdot \rho^{c_1\log(1/\mu_p(f))} \le 0.001 \mu_p(f).
\end{equation}
As $\rho \leq 2\sqrt{p}$, it is sufficient to show that 
\begin{align*}
    2000\frac{1}{\mu_p(f)} \le (\frac{1}{2\sqrt{p}})^{c_1\log(\frac{1}{\mu_p(f)})},
\end{align*}
or equivalently (taking $\log$ from both sides):
\begin{align*}
   \log(2000) + \log(\frac{1}{\mu_p(f)}) \le c_1\log(\frac{1}{\mu_p(f)}) \log(\frac{1}{2\sqrt{p}}).
\end{align*}
This obviously holds, provided that $p$ is smaller than a constant which depends only on $c_1$.

Combining~\eqref{Eq:Product-of-measures1},~\eqref{Eq:Product-of-measures2},~\eqref{Eq:Product-of-measures3}, and~\eqref{Eq:Product-of-measures4}, we get (for all sufficiently small $c_3,p$): 
\begin{align*}
  \mu_p(f)\cdot \mu_\frac{1}{3}(1-g) \le 0.004\mu_p(f).
\end{align*}
This yields $\mu_\frac{1}{3}(1-g) \le 0.004$, or equivalently, $\mu_\frac{1}{3}(g)\ge 0.996$, contradictory to the assumption $\mu_\frac{1}{3}(g)\le 0.99^t$.
\end{proof}

\begin{proof}[Proof of Proposition~\ref{Prop:Global-cross-intersecting-are-small}] 
Let $A,B$ be families that satisfy the assumptions of the proposition, and assume on the contrary that both $\mu_p(A)$ and $\mu_p(B)$ are larger than  $e^{-c_3 \cdot t \cdot \frac{1}{p} \cdot \frac{1}{g^2}}$, where $c_3,p$ are the same as in Lemma~\ref{Lemma:Sharp-threshold-global}. Applying the lemma separately to $A$ and $B$, we get
\[
\mu_{\frac{1}{3}}(A) \ge 0.99^t, \qquad \mbox{and} \qquad \mu_{\frac{1}{3}}(B) \ge 0.99^t.
\]
Since $A$ and $B$ are monotone, it follows from the FKG inequality (specifically, by~\eqref{Eq:FKG} above) that 
\[
\mu_{\frac{1}{3}}(A\cap B) \ge \mu_{\frac{1}{3}}(A)\cdot \mu_{\frac{1}{3}}(B) \ge 0.99^{2t}.
\]
Therefore, $A\cap B$ is a $t$-intersecting family with $\mu_{\frac{1}{3}}(A \cap B) > 0.98^t$, which contradicts the biased Ahlswede-Khachatrian theorem (Theorem~\ref{Thm:Biased-AK} above).
\end{proof}

\section{Large Cross $t$-Intersecting Families Have a Density Bump Inside a Dictatorship}
\label{sec:density-bump}

Informally, a family $\mathcal{F}$ in a product probability space $(K^n,\mu^n)$ is said to \emph{have a density bump} inside a set $S \subset [n]$ if there exists $x \in K^{S}$ s.t. $\mu^{n-|S|}(\mathcal{F}_{S \to x})$ is much larger than $\mu^n(\mathcal{F})$ -- that is, if there exists a restriction which significantly increases the (relative) measure of the family. In these terms, a \emph{global} family is a family that does not admit significant density bumps. 

Several major results in analysis of Boolean functions give sufficient conditions for the existence of a density bump in different settings. (See, for example, the classical Bourgain's theorem~\cite{FriedgutB99} which asserts that any monotone family in $(\{0,1\}^n,\mu_p)$ that has a small total influence, admits a density bump inside a small set $S$, and its recent sharpening by Keevash et al.~\cite{KLLM21}). In this section we show a result of this spirit for large cross $t$-intersecting families.

\begin{prop}\label{Prop:large-cross-intersecting-have-desity-bump}
For any $k \in \mathbb{N}$ and $c>0$, the following holds for all $n \geq n_0(c,k)$ and all $t \leq \frac{n}{10k}$. Let $A,B\subseteq S_n$ be cross $t$-intersecting families. 
If $|A|\cdot|B| > c\cdot(n-t)!^2$, then both $A$ and $B$ have a density bump inside a dictatorship, meaning there exist $i,j$ s.t. $|A_{i\rightarrow j}|>|A|\cdot \frac{k}{n}$, and similarly for $B$. 
\end{prop}

\subsubsection*{Structure of the proof.} The proof goes by contradiction. We assume w.l.o.g. that $B$ does not have a density bump inside a dictatorship and arrive to a contradiction. The proof is divided into three cases, according to the size of $t$ (relative to $n$), and each case consists of four main steps (or part of them):
\begin{enumerate}
    \item \emph{Embedding into a convenient space.} We embed $A,B$ into a more convenient probability space -- usually, $(\{0,1\}^{n^2},\mu_p)$ for a properly chosen value of $p$. In this step, the goal is to show that the embedding does not affect the measure of $A,B$ `too much'.
    
    \item \emph{Replacing one family by a global family.} We replace $A$ by a global restriction $A'$ and show that $B$ has a restriction $B'$ which cross $t$-intersects $A'$ and is not too small. 
    
    \item \emph{Replacing both families by global families.} We replace $B'$ with a global restriction $B''$ and show that $A'$ has a restriction $A''$ which is not too small, is global, and cross $t$-intersects $B''$.
    
    \item \emph{A contradiction to Proposition~\ref{Prop:Global-cross-intersecting-are-small}.} We show that $(A'',B'')$ yield a contradiction to Proposition~\ref{Prop:Global-cross-intersecting-are-small} which asserts that global $t$-cross-intersecting families must be `small'.
\end{enumerate}
As many of the calculations in the three cases are similar, we present each such calculation in detail at the first time it appears and skip the calculation at subsequent appearances. 

\paragraph{\emph{Notation.}} In the course of the proof, we use measures with respect to different spaces. To avoid ambiguity, we mark the measure space in superscript. In particular, $\mu^{S_n}(A)$ and $\mu^{U_n}(A)$ denote the uniform measure of $A$ in $S_n$ and in $U_n=[n]^n$, respectively, while $\mu_p^{N}(\tilde{A})$ and $\mu_p^{S^c}(\tilde{A})$ denote the measure of $\tilde{A}$ in $(\{0,1\}^{n^2},\mu_p)$ and in $(\{0,1\}^{[n^2]\setminus S},\mu_p)$, respectively.

\begin{proof}
Assume on the contrary that for $k,c,n,t$ that satisfy the assumption of the theorem, we have 
\begin{equation}\label{Eq:large-sets}
\mu^{S_n}(A)\cdot \mu^{S_n}(B) \ge c\cdot \frac{(n-t)!^2}{n!^2} \ge c \cdot n^{-2t},
\end{equation}
but $B$ does not have a density bump inside a dictator, meaning that 
\begin{equation}\label{Eq:density-bump}
    \forall i,j: \qquad |B_{i\rightarrow j}| \le |B|\cdot \frac{k}{n}.    
\end{equation}
We divide our proof into three cases:

\begin{itemize}
    \item \emph{Case 1 -- Large $t$}: $\frac{2n}{\log(n)} \leq t \leq n\cdot\frac{1}{10k}$;

    \item \emph{Case 2 -- Medium $t$:} $C(k) \le t \le \frac{2n}{\log(n)}$; 

    \item \emph{Case 3 -- Small $t$}: $t \le \frac{n}{(\log(n))^2}$.
\end{itemize}

\subsection{Case 1 -- Large $t$: $\frac{2n}{\log(n)} \leq t \leq n\cdot\frac{1}{10k}$} 

\paragraph{\emph{Step~1: Embedding.}} We embed $S_n$ into $U_n=[n]^n$ in the standard way (i.e., $\sigma$ goes to $(\sigma(1),\ldots,\sigma(n))$), obtaining two cross $t$-intersecting families $\bar{A},\bar{B} \subset U_n$. By the assumption~\eqref{Eq:large-sets}, we have
\begin{equation}\label{Eq:large-sets11}
    \mu^{U_n}(\bar{A})\cdot \mu^{U_n}(\bar{B}) = \mu^{S_n}(A) \cdot \mu^{S_n}(B) \cdot (\frac{n!}{n^n})^2 \ge c \cdot n^{-2t} e^{-2n}.
\end{equation}
Note that the density bump assumption~\eqref{Eq:density-bump} remains the same in $U_n$, as the embedding affects only the measure of the sets (multiplying it by the same value for all sets) but not their size.

\paragraph{\emph{Step~2: Replacing $\bar{A}$ by a global family.}} Let $\bar{A}'=\bar{A}_{S \to x}$ be a $\sqrt{n}$-global restriction of $\bar{A}$. We claim that $|S|$ is `not too large'. Indeed, by~\eqref{Eq:g-global} and ~\eqref{Eq:large-sets11}, we have 
\begin{equation*}
1 \ge \mu^{([n]\setminus S)^n}(\bar{A}_{S\rightarrow x}) \ge (\sqrt{n})^{|S|}\cdot \mu^{U_n}(A) \ge \sqrt{n}^{|S|} \cdot c n^{-2t}e^{-2n}.    
\end{equation*}
By taking logarithms, we get
\begin{equation*}
\log(1/c) + 2n + 2t\cdot \log(n) \ge |S|\log(\sqrt{n}),    
\end{equation*}
which gives: 
\begin{equation}\label{Eq:large-t-density-bump}
     \frac{|S|}{n} \le \frac{\frac{\log(1/c)}{n} + 2 + \frac{2t\log(n)}{n}}{\frac{1}{2}\cdot \log(n)} \le \frac{1}{2k},
\end{equation}
where the last inequality holds for all $n \geq n_0(c,k)$ since $t\le \frac{n}{10k}$. 

Let $\bar{B}':=\cup_{\{y:y_i \neq x_i ,\forall i \in S\}}(\bar{B}_{S\rightarrow y})$. Clearly, $\bar{B}'$ cross $t$-intersects $\bar{A}'$. We claim that since $|S|$ is not too large, $\bar{B}'$ is non-empty. Indeed, by using a union bound,~\eqref{Eq:density-bump}, and then~\eqref{Eq:large-t-density-bump} we have
\begin{equation*}
|\bar{B}'| \ge |\bar{B}| - \sum_{i \in S} |\bar{B}_{i\rightarrow x_i}| \ge |\bar{B}|\cdot(1-\frac{k}{n}|S|)>0.
\end{equation*}

\paragraph{\emph{Step 3': Reaching a contradiction.}} Let $y \in \bar{B}'$. Every element in $\bar{A}'$ intersects $y$ on at least $t$ coordinates. Hence, $\bar{A}' \subseteq \cup_{|S'| = t}\bar{A}'_{S'\rightarrow y_{S'}}$. By the globalness of $A'$ and a union bound, we have
\begin{equation*}
\begin{split}
\mu^{U_n}(\bar{A}')\cdot n^n &= |\bar{A}'| \le \sum_{|S'|=t}|\bar{A}'_{S'\rightarrow y_{S'}}|=\sum_{|S'|=t}(\mu^{[n]^{[n]\setminus S'}}(\bar{A}'_{S'\rightarrow y_{S'}})\cdot n^{n-t} )\\ 
&\le{n \choose t}\cdot n^{n-t}\sqrt{n}^t\mu^{U_n}(\bar{A}') \le  2^n\cdot n^{n-t}\sqrt{n}^t\mu^{U_n}(\bar{A}').
\end{split}
\end{equation*}
By rearranging and dividing by $\mu^{U_n}(\bar{A}')$ (which is non-zero, since $A$ is non-empty), we get
$2^n \ge (\sqrt{n})^t$, which implies $t < \frac{2n}{\log(n)} $, contradicting our assumption on $t$.

\subsection{Case 2 -- Medium $t$: $C(k) \le t \le \frac{2n}{\log(n)}$}.

\paragraph{\emph{Step~1: Embedding.}} First, we embed $S_n$ into $U_n$ like in Case~1, obtaining $\bar{A},\bar{B} \subset U_n$. As we would like to apply Proposition~\ref{Prop:Global-cross-intersecting-are-small}, we further embed $U_n$ into the space $((\{0,1\}^{n})^n,\mu_p) \sim (\{0,1\}^{n^2},\mu_p)$, by setting $E(i_1,\ldots,i_n)=(e_{i_1},\ldots,e_{i_n})$, where $e_{i_j}$ is the $i_j$'th unit vector in $\{0,1\}^n$. For convenience, we use the notation $N:=n^2$ hereafter.

Note that the embedding preserves the $t$-intersection property. 
In addition, since the image of the embedding is included in the slice $\{x \in \{0,1\}^{N}: \sum x_i = n\}$ on which the $\mu_p$ measure is uniform, the effect of the embedding on the measure of $\bar{A},\bar{B}$ is multiplication by the same factor. We claim that for $p=1/n$, this factor is at least $e^{-n}$. Indeed, for each singleton $i= (i_1,\ldots,i_n) \in U_n$ and its embedding $E(i)=(e_{i_1},\ldots,e_{i_n}) \in \{0,1\}^{N}$, we have  
\begin{equation*}
\mu_p^{N}((e_{i_1},\ldots,e_{i_n})) = p^n\cdot(1-p)^{n^2-n} = \frac{1}{n^n} \cdot (1-\frac{1}{n})^{(n-1)\cdot n} \ge \frac{1}{n^n} \cdot \frac{1}{e^n} = \mu^{U_n}(i) \cdot e^{-n}. \end{equation*}
By additivity of measure, we get $\mu_p^{N}(E(\bar{A})) \geq \mu^{U_n}(\bar{A}) \cdot e^{-n}$ for any family $\bar{A} \subset U_n$.

Let $\hat{A}, \hat{B} \subset \{0,1\}^{N}$ be the embeddings of $\bar{A}, \bar{B}$, respectively. By~\eqref{Eq:large-sets11} we have 
\begin{equation}\label{Eq:large-sets-2}
\mu_p^{N}(\hat{A})\cdot\mu_p^{N}(\hat{B}) \ge e^{-2n} \mu^{U_n}(\bar{A})\mu^{U_n}(\bar{B})
\ge c \cdot e^{-4n}\cdot n^{-2t}.
\end{equation}
Furthermore, the notion of `density bump inside a dictatorship' in the image of $E$ has exactly the same meaning as in~$U_n$, and hence, by~\eqref{Eq:density-bump} we have  
\begin{equation}\label{Eq:density-bump-61}
\frac{k}{n} \cdot \mu_p^{N}(\hat{B}) \ge \mu_p^{N}(\hat{B}_{i\rightarrow 1}).
\end{equation}
Finally, let $\tilde{A},\tilde{B} \subset \{0,1\}^{N}$ be the up-closures of $\hat{A},\hat{B}$ (namely, $\tilde{A} = \{x \in\{0,1\}^{N}: \exists y \in \hat{A}, y \leq x\}$, and similarly for $\tilde{B}$). The families $\tilde{A},\tilde{B}$ are monotone, larger than $\hat{A},\hat{B}$, and cross $t$-intersecting. The family $\tilde{B}$ does not necessarily satisfy the `density bump inside a dictatorship' property, but we will get over that. 

\paragraph{\emph{Step~2: Replacing $\tilde{A}$ by a global family.}}  Let $\tilde{A}'=\tilde{A}_{S \to 1}$ be a $g$-global restriction of $\tilde{A}$, for a value of $g$ that will be determined below. (Note that when we take a global restriction of a monotone family, we can always assume that we restrict the coordinates to $1$'s and not to $0$'s, since restricting a coordinate to $0$ can only decrease the relative measure of the family.)

We claim that $|S|$ is `not too large'. Indeed, by~\eqref{Eq:g-global} and ~\eqref{Eq:large-sets-2}, we have 
\begin{equation*}
1 \ge \mu_p^{S^c}(\tilde{A}') \ge g^{|S|}\cdot \mu_p^{N}(\tilde{A}) \ge g^{|S|} \cdot c n^{-2t}e^{-4n}.    
\end{equation*}
By the same computation as in~\eqref{Eq:large-t-density-bump}, using the assumption $t \leq \frac{2n}{\log(n)}$, this yields 
\begin{equation}\label{Eq:medium-t-density-bump}
|S| \leq \frac{n}{2k},    
\end{equation}
provided that $g \geq g_0(k)$ and $n \geq n_0(c,k)$. 

Let $\tilde{B}':=\tilde{B}_{S\rightarrow 0} \supseteq \hat{B}_{S\rightarrow 0}$. We claim that since $|S|$ is `not too large', $\tilde{B}'$ is `not too small'. Indeed, by a union bound,~\eqref{Eq:density-bump-61}, and~\eqref{Eq:medium-t-density-bump}, we have
\begin{equation*}
\mu_p^{N}(\tilde{B}') \ge \mu_p^{N}(\hat{B}_{S\rightarrow 0}) > \mu_p^{N}(\hat{B}) - \sum_{i \in S} \mu_p^{N}(\hat{B}_{i\rightarrow 1}) \ge \mu_p^{N}(\hat{B})\cdot(1-\frac{k}{n}|S|) 
\ge \frac{1}{2} \mu_p^{N}(\hat{B}).
\end{equation*}
As $\mu_p^{S^c}(\tilde{B}') > \mu_p^{N}(\tilde{B}')$,~\eqref{Eq:large-sets-2} yields
\begin{equation}\label{Eq:large-sets-3}
\mu_p^{S^c}(\tilde{A}')\cdot \mu_p^{S^c}(\tilde{B}') \geq \frac{g^{|S|}}{2}\mu_p^{N}(\hat{A})\cdot \mu_p^{N}(\hat{B}) \ge \frac{c}{2}\cdot e^{-4n}\cdot n^{-2t}.
\end{equation} 
Thus, the families  $\tilde{A}',\tilde{B}'$ we obtain are monotone, cross $t$- intersecting, `not too small', and one of them (i.e., $\tilde{A}'$) is global. 

\paragraph{\emph{Step~3: Replacing $\tilde{A}$ and $\tilde{B}$ by global families.}} Let $\tilde{B}"=\tilde{B}'_{S' \to 1}$ be a $g'$-global restriction of $\tilde{B}'$, for $g'$ that will be determined below. By~\eqref{Eq:g-global} and ~\eqref{Eq:large-sets-3}, we have 
\begin{equation*}
1 \ge \mu_p^{(S \cup S')^c}(\tilde{B}'') \ge (g')^{|S'|}\cdot \mu_p^{S^c}(\tilde{B}') \ge (g')^{|S'|} \cdot \frac{c}{2}\cdot  n^{-2t}e^{-4n}.    
\end{equation*}
By the same computation as in~\eqref{Eq:medium-t-density-bump}, this yields 
\begin{equation}\label{Eq:medium-t-density-bump2}
|S| \leq \frac{n}{2g},    
\end{equation}
provided that $g' \geq g_1(g)$ and $n \geq n_0(c,k)$. 

Let $\tilde{A}'':=\tilde{A}'_{S'\rightarrow 0}$. By a union bound, the globalness of $\tilde{A}'$, and~\eqref{Eq:medium-t-density-bump2}, we have
\begin{equation}\label{Eq:dens76}
\mu_p^{S^c}(\tilde{A}") \ge \mu_p^{S^c}(\tilde{A}') - \sum_{i \in S} \mu_p^{S^c}(\tilde{A}'_{i\rightarrow 1}) 
\ge \mu_p^{S^c}(\tilde{A}')\cdot(1-g\cdot p\cdot|S'|)
\ge \frac{1}{2}\mu_p^{S^c}(\tilde{A}').
\end{equation}
As $\mu_p^{(S \cup S')^c}(\tilde{A}'') > \mu_p^{S^c}(\tilde{A}'')$,~\eqref{Eq:large-sets-3} yields
\begin{equation}\label{Eq:large-sets-4}
\mu_p^{(S \cup S')^c}(\tilde{A}'')\cdot \mu_p^{(S \cup S')^c}(\tilde{B}'') \geq \frac{(g')^{|S'|}}{2}\mu_p^{S^c}(\tilde{A}')\cdot \mu_p^{S^c}(\tilde{B}')
\ge \frac{c}{4}\cdot e^{-4n}\cdot n^{-2t},
\end{equation} 
meaning that $\tilde{A}'',\tilde{B}''$ are `not too small'. We claim that $\tilde{A}" \subset \{0,1\}^{[N] \setminus (S \cup S')}$ is $2g$-global. As $\tilde{A}''$ is monotone, we only have to show that for any $T \subset [n^2] \setminus (S\cup S')$ we have
\[
\mu_p^{(S\cup S'\cup T)^c}(\tilde{A}"_{T\rightarrow 1}) \le (2g)^{|T|}\mu_p^{(S \cup S')^c}(\tilde{A}'').
\]
This indeed holds, as
\begin{equation}\label{Eq:2g-global}
\begin{split}
\mu_p^{(S\cup S'\cup T)^c}(\tilde{A}"_{T\rightarrow 1}) &= \mu_p^{(S\cup S'\cup T)^c}(\tilde{A}'_{S' \rightarrow 0, T\rightarrow 1})  \le \mu_p^{(S\cup T)^c}(\tilde{A}'_{T\rightarrow 1}) \le g^{|T|}\mu_p^{S^c}(\tilde{A}') \\&\le g^{|T|}\cdot 2\mu_p^{S^c}(\tilde{A}") \le (2g)^{|T|}\mu_p^{(S\cup S')^c}(\tilde{A}")
\end{split}
\end{equation}
where in the first inequality we used the monotonicity of $\tilde{A}'$, in the second one we used its globalness, and in the third one we used~\eqref{Eq:dens76}.
Hence, the families  $\tilde{A}'',\tilde{B}''$ we obtain are monotone, cross $t$- intersecting, `not too small', and $G$-global, for $G=\max(2g,g')$. 

\paragraph{\emph{Step~4: Reaching a contradiction by applying Proposition~\ref{Prop:Global-cross-intersecting-are-small} to $\tilde{A}'',\tilde{B}''$.}} 
By~\eqref{Eq:large-sets-4}, we have
\[
\mu_p^{(S\cup S')^c}(\tilde{A}"),\mu_p^{(S\cup S')^c}(\tilde{B}") \ge \frac{c}{4}\cdot e^{-4n}\cdot n^{-2t} \ge  e^{-c_3 \cdot t \cdot \frac{1}{\frac{1}{n}} \cdot \frac{1}{G^2}},
\]
where the second inequality holds for all $t>C(g,c_3)$ and $n \ge n_0(c,c_3,k)$, by a calculation similar to~\eqref{Eq:large-t-density-bump}. (Note that as $g$ depends only on $k$, the condition on $t$ boils down to $t >C(k)$.) Hence, $(\tilde{A}'',\tilde{B}'')$ contradict the assertion of Proposition~\ref{Prop:Global-cross-intersecting-are-small}, being a pair of large monotone cross $t$-intersecting families in $(\{0,1\}^{N \setminus (S\cup S')},\mu_{1/n})$. 

\subsection{Case 3 -- Small $t$: $t \le \frac{n}{(\log(n))^2}$}. The complex part in this case is the embedding step, presented in detail below. The other steps are similar to the corresponding steps of Case~2 and will be presented briefly.

\paragraph{\emph{Step~1: Embedding.}} Unlike Cases~1,2, we do not embed $S_n$ into $U_n$, since this results in a loss of a factor of $e^{-2n}$ in the measure, which we cannot afford in this case. Instead, we embed $S_n$ directly into the space $((\{0,1\}^{n})^n,\mu_p) \sim (\{0,1\}^{n^2},\mu_p)$, by setting $E(\sigma)=(e_{\sigma(1)},\ldots,e_{\sigma(n)})$, where $e_{i}$ is the $i$'th unit vector in $\{0,1\}^n$. Note that the embedding preserves the $t$-intersection property. 

Let $\hat{A}, \hat{B} \subset \{0,1\}^{N}$ be the embeddings of $A, B$, respectively, and let $\tilde{A},\tilde{B} \subset \{0,1\}^{N}$ be the up-closures of $\hat{A},\hat{B}$ (namely, $\tilde{A} = \{x \in\{0,1\}^{N}: \exists y \in \hat{A}, y \leq x\}$, and similarly for $\tilde{B}$). When convenient, we view $\hat{A},\hat{B},\tilde{A},$ and $\tilde{B}$ as families of subsets of $[N]$, using the natural correspondence between elements of $\{0,1\}^N$ and subsets of $[N]$. The families $\tilde{A},\tilde{B}$ are monotone
and cross $t$-intersecting. We claim that their biased measure is not significantly less than the measure of $A,B$ (respectively) in $S_n$.
\begin{claim}
In the above definitions, for any $p\ge \frac{10\log(n)}{n}$ we have
\begin{equation}\label{Eq:coupling-measure}
\mu_p^N(\tilde{A})\geq \frac{1}{2}\mu^{S_n}(A), \qquad \mbox{and} \qquad \mu_p^N(\tilde{B})\geq \frac{1}{2}\mu^{S_n}(B).
\end{equation}  
\end{claim}

\begin{proof}[Proof of the claim]
We define a coupling between $(S_n,\mu)$ and $(\{0,1\}^N,\mu_p)$, in the following way. We draw $x \in \{0,1\}^N$ according to $\mu_p$, and choose $\sigma \in S_n$ such that $E(\sigma) \leq x$ coordinate-wise, at uniform distribution from all the prospects. If there are no prospects, we choose $\sigma$ uniformly from $S_n$.
By symmetry, this is indeed a coupling between $(S_n,\mu)$ and $(\{0,1\}^N,\mu_p)$.

Consider two indicator variables on the coupling probability space:
\[
X(x,\sigma) = 1\{x \in \tilde{A}\}, \qquad Y(x,\sigma)= 1\{\sigma \in A\}.
\]
Let $U \subset \{0,1\}^n$ be the up-closure of $E(S_n)$. Namely, $U = \{x \in\{0,1\}^{N}: \exists \sigma \in S_n, E(\sigma) \leq x\}$. Note that by symmetry, the distribution of $(\sigma|x \in U)$ is the uniform distribution on $S_n$. 

Observe that
\begin{align*}
    \mu^{S_n}(A) = \mathbb{E}(Y|x \in U) \le  \mathbb{E}(X|x \in U), 
\end{align*}
where the inequality holds since when we condition on $x \in U$, we have $Y\le X$ as $Y=1$ implies $X=1$. On the other hand,  
\begin{align*}
    \mu_p^N(\tilde{A}) = \mathbb{E}(X) \ge \mathbb{E}(X|x \in U) \cdot \mu_p^N(U).   
\end{align*}
Hence,~\eqref{Eq:coupling-measure} will follow once we show that 
\begin{equation}\label{Eq:Closure-of-Sn}
    \mu_p^N(U) \ge \frac{1}{2}.
\end{equation}
To see this, consider the natural division of the $n^2$ coordinates of $\{0,1\}^N$ into $n$ consecutive $n$-tuples $T_1,\ldots,T_n$. For each $y=(y^1,y^2,\ldots,y^n) \in \{0,1\}^N$, we may define a bipartite graph $G_y \subset K_{n,n}$ in which the vertices of the left side correspond to $T_1,\ldots,T_n$, the vertices of the right side correspond to $1,2,\ldots,n$, and for each $i$, the vertex that corresponds to $T_i$ is connected to all vertices $\{j \in [n]:y^i_j=1\}$. 

Observe that if $y \not \in U$, then $G_y$ does not contain a perfect matching. By Hall's marriage theorem, this implies that for some $1 \leq k \leq n$, there exist $k$ subsets $T_{i_1},\ldots,T_{i_k}$ such that in $G_y$, the vertices of $\cup_{j=1}^k T_{i_j}$ have at most $k-1$ neighbors in total. Taking the minimal $k$ with this property, we may assume that they have exactly $k-1$ neighbors in total. The probability of this event (where $y \sim (\{0,1\}^N,\mu_p$) is bounded from above by
\begin{equation*}
{n \choose {k-1}} \cdot(1-p)^{k\cdot(n-k+1)}.
\end{equation*}
Hence, by a union bound we have
\begin{equation*}
\begin{split}
    \mu_p(U^c) &\leq 
    \sum_{k=1}^n {n \choose k}{n \choose {k-1}} \cdot(1-p)^{k\cdot(n-k+1)} \\
    &\le 2\sum_{k=1}^{\lceil \frac{n+1}{2} \rceil} {n \choose k}{n \choose {k-1}} \cdot(1-\frac{10\log(n)}{n})^{k\cdot\lceil \frac{n+1}{2} \rceil} \le 2\sum_{k=1}^{\frac{n+1}{2}} n^{2k+1}\cdot (\frac{1}{n^{5}})^k \le \frac{1}{2}.
\end{split}
\end{equation*}
This proves~\eqref{Eq:Closure-of-Sn}, and thus, completes the proof of the claim.
\end{proof}

Returning to the embedding step, we fix $p= \frac{10\log(n)}{n}$, embed $A$ and $B$ into $(\{0,1\}^N,\mu_p)$ to obtain $\hat{A},\hat{B}$, and then we consider the up-closures $\tilde{A},\tilde{B}$ of $\hat{A},\hat{B}$. As was mentioned above, $\tilde{A},\tilde{B} \subset \{0,1\}^N$ are monotone and cross $t$-intersecting. By~\eqref{Eq:large-sets} and~\eqref{Eq:coupling-measure}, we have
\begin{equation}\label{Eq:large-sets-32}
\mu_p^{N}(\tilde{A})\cdot \mu_p^{N}(\tilde{B}) > \frac{c}{4} n^{-2t}. 
\end{equation}


\paragraph{\emph{Step~2: Replacing $\tilde{A}$ by a global family.}}  Let $\tilde{A}'=\tilde{A}_{S \to 1}$ be a $g$-global restriction of $\tilde{A}$, for a value of $g$ that will be determined below. 

We claim that $|S|$ is `not too large'. Indeed, by the same computation as in the corresponding step of Case~2 (specifically,~\eqref{Eq:medium-t-density-bump}), we obtain 
\begin{equation}\label{Eq:small-t-density-bump}
|S| \leq \frac{n}{2k},    
\end{equation}
provided that $g \geq g_0(k)$ and $n \geq n_0(c,k)$.

Let $\tilde{B}':=\tilde{B}_{S\rightarrow 0}$. We claim that since $|S|$ is `not too large', $\tilde{B}'$ is `not too small'. To see this, let $\hat{B}' := \hat{B}_{S\rightarrow 0}$, where by setting the $((i-1)\cdot n+j)$'th coordinate to $0$ we mean that we look at the embeddings of all permutations $\sigma$ s.t.~$\sigma(i) \neq j$. By a union bound and~\eqref{Eq:small-t-density-bump}, we have 
\begin{equation}\label{Eq:density-bump-31}
|\hat{B'}| \ge |\hat{B}| - \sum_{i \in S} |\hat{B}_{i\rightarrow 0}| \ge |\hat{B}|\cdot(1-\frac{k}{n}|S|) \ge \frac{1}{2}|\hat{B}|.   
\end{equation}
(Note that $|\hat{B}|$ is the number of elements of $\hat{B'}$, viewed as a family of subsets of $[N]$, and similarly for the other set sizes here.) Let $R'\subset \{0,1\}^N$ be the up-closure of $\hat{B'}$. As $R'_{S\rightarrow 0} = \tilde{B}'$, we have
\begin{align*}
    \mu^{S_n}(B) = \frac{|\hat{B}|}{n!} \le 2\frac{|\hat{B'}|}{n!} \le 4\mu_p^{N}(R') = 4\mu_p^{S^c}(R'_{S\rightarrow 0 }) = 4\mu_p^{S^c}(\tilde{B}'),
\end{align*}
where the first inequality is by~\eqref{Eq:density-bump-31}, the second inequality is by~\eqref{Eq:coupling-measure}, and the third inequality holds since the coordinates of $S$ are off in $\hat{B'}$.

Finally,~\eqref{Eq:large-sets} and~\eqref{Eq:coupling-measure} yield
\begin{equation}\label{Eq:large-sets-33}
\mu_p^{S^c}(\tilde{A}')\cdot \mu_p^{S^c}(\tilde{B}') \geq \frac{g^{|S|}}{4}\mu_p^{N}(\tilde{A})\cdot \mu^{S_n}(B) \ge \frac{c}{8}\cdot n^{-2t}.
\end{equation} 
Thus, the families  $\tilde{A}',\tilde{B}'$ we obtain are monotone, cross $t$- intersecting, `not too small', and one of them (i.e., $\tilde{A}'$) is global.

\paragraph{\emph{Step~3: Replacing $\tilde{A}$ and $\tilde{B}$ by global families.}} Let $\tilde{B}"=\tilde{B}'_{S' \to 1}$ be a $g'$-global restriction of $\tilde{B}'$, for $g'$ that will be determined below. By~\eqref{Eq:g-global} and ~\eqref{Eq:large-sets-33}, we have 
\begin{equation*}
1 \ge \mu_p^{(S \cup S')^c}(\tilde{A}'') \ge (g')^{|S'|}\cdot \mu_p^{S^c}(\tilde{A}') \ge (g')^{|S'|} \cdot \frac{c}{8}\cdot  n^{-2t}.    
\end{equation*}
By the same computation as in~\eqref{Eq:medium-t-density-bump}, using the assumption $t\leq \frac{n}{\log(n)^2}$, this yields 
\begin{equation}\label{Eq:small-t-density-bump2}
|S| \leq \frac{1}{p \cdot 2g}    
\end{equation}
provided that $g' \geq g_1(g)$ and $n \geq n_0(c,k)$. 

Let $\tilde{A}'':=\tilde{A}'_{S'\rightarrow 0}$. By a union bound, the globalness of $\tilde{A}'$, and~\eqref{Eq:small-t-density-bump2}, we have
\begin{equation}\label{Eq:dens761}
\mu_p^{S^c}(\tilde{A}") \ge \mu_p^{S^c}(\tilde{A}') - \sum_{i \in S} \mu_p^{S^c}(\tilde{A}'_{i\rightarrow 1}) 
\ge \mu_p^{S^c}(\tilde{A}')\cdot(1-g\cdot p\cdot|S'|)
\ge \frac{1}{2}\mu_p^{S^c}(\tilde{A}').
\end{equation}
As $\mu_p^{(S \cup S')^c}(\tilde{A}'') > \mu_p^{S^c}(\tilde{A}'')$,~\eqref{Eq:large-sets-33} yields
\begin{equation}\label{Eq:large-sets-34}
\mu_p^{(S \cup S')^c}(\tilde{A}'')\cdot \mu_p^{(S \cup S')^c}(\tilde{B}'') \geq \frac{(g')^{|S'|}}{2}\mu_p^{S^c}(\tilde{A}')\cdot \mu_p^{S^c}(\tilde{B}')
\ge \frac{c}{16} \cdot n^{-2t},
\end{equation} 
meaning that $\tilde{A}'',\tilde{B}''$ are `not too small'. In addition, by exactly the same argument as in Case~2 (specifically, by~\eqref{Eq:2g-global}), $\tilde{A}" \subset \{0,1\}^{[N] \setminus (S \cup S')}$ is $2g$-global. Hence, the families  $\tilde{A}'',\tilde{B}''$ we obtain are monotone, cross $t$-intersecting, `not too small', and $G$-global, for $G=\max(2g,g')$.

\paragraph{\emph{Step~4: Reaching a contradiction by applying Proposition~\ref{Prop:Global-cross-intersecting-are-small} to $\tilde{A}'',\tilde{B}''$.}} 
By~\eqref{Eq:large-sets-34}, we have
\[
\mu_p^{(S\cup S')^c}(\tilde{A}"),\mu_p^{(S\cup S')^c}(\tilde{B}") \ge \frac{c}{16}\cdot n^{-2t} \ge  e^{-c_3 \cdot t \cdot \frac{1}{\frac{10\log(n)}{n}} \cdot \frac{1}{G^2}},
\]
where the second inequality holds for all $t \geq 1$ and $n \ge n_0(c,c_3,k)$, by a calculation similar to~\eqref{Eq:large-t-density-bump}. Hence, $(\tilde{A}'',\tilde{B}'')$ contradict the assertion of Proposition~\ref{Prop:Global-cross-intersecting-are-small}, being a pair of large monotone cross $t$-intersecting families in $(\{0,1\}^{N \setminus (S\cup S')},\mu_{10\log(n)/n})$. 
\end{proof}

\section{Proof of the Main Theorem}
\label{sec:main}

In this section we present the proof of Theorem~\ref{Thm:Main}. We prove a more general statement, which can be proved more conveniently by induction. 
\begin{prop} \label{Prop:main}
There exists a constant $c_0>0$ such that the following holds for any $t \in \mathbb{N}$.
\begin{enumerate}
    \item For any $n \geq \lfloor c_0 \cdot t \rfloor$, if $A,B\subseteq S_n$ are cross $t$-intersecting families, then
    \[
    |A||B| \le (n-t)!^2 \cdot \max\{4^{2\lfloor c_0 \cdot t \rfloor-n},1\};
    \]
    \item For any $n \geq 2\lfloor c_0 \cdot t \rfloor$, if $A,B\subseteq S_n$ are cross $t$-intersecting families and $|A||B| \ge \frac{3}{4}(n-t)!^2$, then there exist $i_1,i_2,\ldots,i_t$ and $j_1,j_2,\ldots,j_t$ such that $A,B \subset (S_n)_{i_1\to j_1,i_2\to j_2,\ldots, i_t\to j_t}$. If in addition, $A=B$, then the same assertion holds under the weaker assumption $|A| \ge \frac{3}{4}(n-t)!$ 
\end{enumerate}
\end{prop}
Assertion~(2) of Proposition~\ref{Prop:main} implies Theorem~\ref{Thm:Main}. 

\begin{rem}\label{Rem:Optimal-constant}
We note, as will be apparent from the proof, that the constant $\frac{3}{4}$ in the `stability' statement (both for general $A,B$ and for $A=B$) can be replaced by any constant larger than $1-1/e$, provided that the constant $c_0$ is adjusted accordingly.
\end{rem}

\begin{proof}
We prove that the proposition holds with 
\begin{equation}\label{Cond:c0}
    c_0 \geq \max\{n_0(\frac{2}{3},50),500\}, 
\end{equation}
where $n_0$ is as defined in the statement of Proposition~\ref{Prop:large-cross-intersecting-have-desity-bump}. We prove Assertion~(1) of the proposition by induction on $t$ and $n$. In the induction step we assume that Assertion~(1) holds for all $(t',n')$ such that either $(t'<t) \wedge (\lfloor c_0 t' \rfloor \leq n')$ or $(t'=t)\wedge (\lfloor c_0 t \rfloor \leq n'<n)$ and prove Assertion~(1) for $(t,n)$. Assertion~(2) will follow easily from the proof of Assertion~(1), as we will show at the end of the proof.

\medskip

\paragraph{\emph{Induction basis.}} In the basis case, we have to show that for any $t \in \mathbb{N}$, Assertion~(1) holds for $n=\lfloor c_0 t \rfloor$. That is, we have to show that if $n=\lfloor c_0 t \rfloor$ and $A,B \subset S_n$ are cross $t$-intersecting families, then $|A||B| \le 4^{n}(n-t)!^2$. 

We may assume that $A,B \neq \emptyset$. Let $f \in A$. Every element $g\in B$ agrees with $f$ on at least $t$ coordinates (when we view them as elements of $[n]^n$). For each choice of $t$ coordinates, there are at most $(n-t)!$ elements of $S_n$ which agree with $f$ on those $t$ coordinates. Hence,
\[
|B| \le {n \choose t}\cdot (n-t)! \le 2^n(n-t)!. 
\]
By the same argument, we have $|A| \le 2^n(n-t)!,$ and thus, $|A||B| \le 4^{n}(n-t)!^2$, as asserted. 

\medskip 

\paragraph{\emph{Induction step.}} The proof of the induction step consists of two steps:
\begin{enumerate}
    \item \emph{Both families are almost contained in the same dictatorship.} We show that assuming the induction hypothesis, we can `upgrade' Proposition~\ref{Prop:large-cross-intersecting-have-desity-bump} to show that both $A$ and $B$ are almost fully contained in the same dictatorship $(S_n)_{i\to j}$.

    \item \emph{Bootstrapping.} We show that the first step can be applied repeatedly to show that both $A$ and $B$ are almost fully contained in the same $t$-umvirate $(S_n)_{i_1\to j_1,\ldots,i_t \to j_t}$, and deduce Assertion~(1) for $(t,n)$ from this structural result.
\end{enumerate}

\subsection*{Step~1: Both families are almost contained in the same dictatorship}

We prove the following.
\begin{claim}\label{Claim:Almost-in-dictatorship}
    Let $c_0 \geq \max\{n_0(\frac{2}{3},50),500\}$, where $n_0$ is as defined in the statement of Proposition~\ref{Prop:large-cross-intersecting-have-desity-bump}. Let $t \in \mathbb{N}$ and let $n > \lfloor c_0 t \rfloor$. Assume that Assertion~(1) of Proposition~\ref{Prop:main} holds for all $(t',n')$ such that either $(t'<t) \wedge (\lfloor c_0 t' \rfloor \leq n')$ or $(t'=t)\wedge (\lfloor c_0 t \rfloor \leq n'<n)$. 
    Let $A,B \subset S_n$ be cross $t$-intersecting families such that   
    \begin{equation}\label{Eq:large-sets-41}
    |A||B| > r(n,t)\cdot\frac{2}{3}\cdot (n-t)!^2,    
    \end{equation}
    where $r(n,t) =\max\{4^{2\lfloor c_0 t \rfloor-n},1\}$.
    Then there exist $i,j \in [n]$ such that 
    \begin{equation}\label{Eq:Almost-in-dictatorship}
        |A_{i \rightarrow j}| \ge |A|(1-\frac{7}{n}), \qquad \mbox{and} \qquad |B_{i \rightarrow j}| \ge |B|(1-\frac{7}{n}).
    \end{equation}
\end{claim}

\begin{proof}[Proof of the claim]
    By Proposition~\ref{Prop:large-cross-intersecting-have-desity-bump}, applied with $c=\frac{2}{3}$ and $k=50$, $A$ has a density bump into a dictator, meaning w.l.o.g.~that 
\begin{equation}\label{Eq:a-ge-k}
|A_{1\rightarrow i}| = \frac{a\cdot|A|}{n}, \qquad \mbox{for some} \qquad a>50.
\end{equation}
Note that Proposition~\ref{Prop:large-cross-intersecting-have-desity-bump} can be indeed applied to $A$ and $B$, since the condition on $c_0$ guarantees that the hypotheses $n \geq n_0(c,k)$ and $t \leq \frac{n}{10k}$ of the proposition are satisfied.

For each $j \neq i$, the families $A_{1\rightarrow i}, B_{1\rightarrow j} \subset [n]^n$ are cross $t$-intersecting. We can view both families as residing in `copies' of $S_{n-1}$ by viewing the restrictions of all elements of $A_{1\rightarrow i}, B_{1\rightarrow j}$ to the domain $\{2,3,\ldots,n\}$. (Note that the set of bijections between $\{2,3,\ldots,n\}$ and $[n] \setminus \{i\}$ is clearly isomorphic to $S_{n-1}$). Furthermore, we can view those `copies' of $S_{n-1}$ as the same copy, by renaming the element $j$ in the range of $A_{1\to i}$ to $i$ (which does not affect the $t$-intersection property). Hence, we can apply the induction hypothesis to $A_{1\rightarrow i},B_{1\rightarrow j}$ to deduce
\[
|A_{1\rightarrow i}||B_{1\rightarrow j}| \leq r(n-1,t)\cdot (n-t-1)!^2.
\]
Using a union bound and noting that $r(n-1,t) \le 4r(n,t)$, we get
\begin{equation}\label{Eq:in}
\frac{a}{n}|A||\cup_{j \neq i}(B_{1\rightarrow j})|=|A_{1\rightarrow i}||\cup_{j \neq i}(B_{1\rightarrow j})| \le (n-1)\cdot 4r(n,t)\cdot (n-t-1)!^2.
\end{equation}
Equation~\eqref{Eq:in} allows obtaining a lower bound on $|B_{1\to i}|$. Indeed, by combining~\eqref{Eq:in} with the assumption~\eqref{Eq:large-sets-41}, we have 
\begin{align*}
\begin{split}
    \frac{2}{3}\cdot r(n,t)(n-t)!^2 &\leq |A|\cdot|B| = |A|\cdot (|\cup_{j \neq i}(B_{1\rightarrow j})|+|B_{1\rightarrow i}|) \\
    &\leq |A|\cdot|B_{1\rightarrow i}|+\frac{n}{a}\cdot4r(n,t)(n-1)\cdot(n-t-1)!^2,
\end{split}
\end{align*}
and consequently,
\begin{equation}\label{Eq:B-density}
    |B_{1\rightarrow i}| \ge \frac{r(n,t)(n-t-1)!^2\cdot[\frac{2}{3}\cdot(n-t)^2-\frac{4n(n-1)}{a}]}{|A|}.
\end{equation}
Now, we obtain an upper bound on $|\cup_{j \neq i}(A_{1\rightarrow j})|$. By applying the argument we used to obtain~\eqref{Eq:in} with the roles of $A,B$ reversed, we have $|B_{1\rightarrow i}||\cup_{j \neq i}(A_{1\rightarrow j})| \le (n-1)\cdot 4r(n,t) \cdot (n-t-1)!^2$. Combining this with~\eqref{Eq:B-density}, we get:
\begin{align*}
    |\cup_{j \neq i}(A_{1\rightarrow j})| \le \frac{(n-1)\cdot 4r(n,t)\cdot (n-t-1)!^2\cdot|A|}{r(n,t)(n-t-1)!^2\cdot[\frac{2}{3}\cdot(n-t)^2-\frac{4n(n-1)}{a}]},
\end{align*}
and consequently, 
\begin{equation}\label{Eq:a0}
\frac{|\cup_{j \neq i}(A_{1\rightarrow j})|}{|A|} \le \frac{4(n-1)}{\frac{2}{3}\cdot(n-t)^2-\frac{4n(n-1)}{a}}.    
\end{equation}
Equation~\eqref{Eq:a0} allows us obtaining an inequality involving only $t,n,$ and $a$:
\begin{equation}\label{Eq:a}
    \frac{a}{n} = \frac{|A_{1\rightarrow i}|}{|A|} = 1 -\frac{|\cup_{j \neq i}(A_{1\rightarrow j})|}{|A|} \ge 1 - \frac{4(n-1)}{\frac{2}{3}\cdot(n-t)^2-\frac{4n(n-1)}{a}}.
\end{equation}
Writing~\eqref{Eq:a} as a quadratic inequality in the variable $a$ (viewing $n$ and $t$ as parameters) and solving it, we obtain
\begin{equation*}
a \ge \frac{n}{2} \cdot (1+\sqrt{\Delta}) \qquad \mbox{or} \qquad a \le \frac{n}{2} \cdot (1-\sqrt{\Delta}), 
\end{equation*}
where $\Delta = 1-24\frac{n-1}{(n-t)^2}$. We claim that the latter is impossible. Indeed, as we assume $c_0 \geq 500$ (see~\eqref{Cond:c0}), we have $n \geq c_0 t \geq 500t \geq 500$, and thus, 
\begin{equation}\label{Eq:a3}
\sqrt{\Delta}=\sqrt{1-24\frac{n-1}{(n-t)^2}} \geq \sqrt{1-24\frac{n}{(499n/500)^2}} \geq \sqrt{1-\frac{25}{n}}>1-\frac{100}{n}.
\end{equation}
Hence, $a \le \frac{n}{2} \cdot (1-\sqrt{\Delta})$ would imply
\[
a \le \frac{n}{2} \cdot (1-\sqrt{\Delta}) < \frac{n}{2}(1-(1-\frac{100}{n}))=50,
\]
in contradiction to~\eqref{Eq:a-ge-k}. Therefore, we have
\[
a \ge \frac{n}{2} \cdot (1+\sqrt{\Delta}) > \frac{n}{2}(1+(1-\frac{100}{n}))=n-50,
\]
or equivalently,
\begin{equation}\label{Eq:a4}
|A_{1\rightarrow i}|\ge (1-\frac{50}{n})|A|.
\end{equation}
To see that an even stronger statement holds for $B_{1\to i}$, we combine~\eqref{Eq:large-sets-41},~\eqref{Eq:in}, and~\eqref{Eq:a4} to get
\begin{align*}
\begin{split}
    1 - \frac{|B_{1\rightarrow i}|}{|B|} &= \frac{|\cup_{j \neq i}(B_{1\rightarrow j})|}{|B|} = \frac{|A_{1\to i}| |\cup_{j \neq i}(B_{1\rightarrow j})|}{|A_{1\to i}| |B|} \\
    &\le \frac{(n-1)\cdot 4r(n,t)(n-t-1)!^2}{ (1-\frac{50}{n})r(n,t)\frac{2}{3}(n-t)!^2} =\frac{6n(n-1)}{(n-50)(n-t)^2} \\
    &\le \frac{6n^2}{(n-n/10)(n-n/500)^2} \le \frac{7}{n},
\end{split}
\end{align*}
where the penultimate inequality holds since $n \ge 500 t \geq 500$ by~\eqref{Cond:c0}. Consequently, we have $|B_{1\rightarrow i}|\ge (1-\frac{7}{n})|B|$, and by the same argument with the roles of $A,B$ reversed, we have $|A_{1\rightarrow i}|\ge (1-\frac{7}{n})|A|$ as well. This completes the proof of the claim.
\end{proof}

\subsection*{Step~2: Bootstrapping} Recall that we assume that Assertion~(1) of Proposition~\ref{Prop:main} holds for all $(t',n')$ such that either $(t'<t) \wedge (\lfloor c_0 t' \rfloor \leq n')$ or $(t'=t)\wedge (\lfloor c_0 t \rfloor \leq n'<n)$, and we want to prove Assertion~(1) for $(t,n)$. 

Let $A,B \subset S_n$ be cross $t$-intersecting families such that $|A||B| \geq r(n,t)(n-t)!^2$. Assuming that $c_0$ satisfies~\eqref{Cond:c0}, we may apply Claim~\ref{Claim:Almost-in-dictatorship} to deduce that there exist $i_1,j_1 \in [n]$ such that 
\begin{equation}\label{Eq:Boot1}
|A_{i_1 \rightarrow j_1}| \ge |A|(1-\frac{7}{n}), \qquad \mbox{and} \qquad |B_{i_1 \rightarrow j_1}| \ge |B|(1-\frac{7}{n}).
\end{equation}
Due to the assumption $|A||B| \geq r(n,t)(n-t)!^2$, this implies
\[
|A_{i_1\to j_1}||B_{i_1\to j_1}| \geq (1-\frac{7}{n})^2 r(n,t)(n-t)!^2 > \frac{2}{3} \cdot r(n,t)(n-t)!^2.
\]  
The families $A_{i_1\rightarrow j_1}, B_{i_1\rightarrow j_1}$ are cross $(t-1)$-intersecting, and may be viewed as subfamilies of the same copy of $S_{n-1}$. Hence, we may apply Claim~\ref{Claim:Almost-in-dictatorship} to these restrictions, to deduce that there exist $i_2,j_2 \in [n]$ such that  
\begin{equation*}
|A_{i_1 \rightarrow j_1,i_2\to j_2}| \ge |A_{i_1\to j_1}|(1-\frac{7}{n}), \mbox{ and } |B_{i_1 \rightarrow j_1,i_2\to j_2}| \ge |B_{i_1\to j_1}|(1-\frac{7}{n})
\end{equation*}
and
\[
|A_{i_1\to j_1,i_2\to j_2}||B_{i_1\to j_1,i_2 \to j_2}| \geq (1-\frac{7}{n})^4 r(n,t)(n-t)!^2 > \frac{2}{3} \cdot r(n,t)(n-t)!^2.
\] 
(Note that the claim can indeed be applied again, as $r(n-1,t-1) \leq r(n,t)$, $((n-1)-(t-1))!=(n-t)!$, and the inductive hypothesis for $(t-1,n-1)$ is included in the inductive hypothesis for $(t,n)$). Since by~\eqref{Cond:c0}, 
\begin{equation}\label{Eq:Boot1.5}
(1-\frac{7}{n})^{2t}\geq (1-\frac{7}{500t})^{2t} \geq (e^{-\frac{14}{500t}})^{2t} > 0.94 > \frac{2}{3},
\end{equation}
we may apply Claim~\ref{Claim:Almost-in-dictatorship} again and again $t$ times, to deduce that there exist $i_1,\ldots,i_t$ and $j_1,\ldots,j_t$ such that for any $1 \leq \ell \leq t-1$,
\begin{equation}\label{Eq:Boot2}
|A_{i_1 \rightarrow j_1,\ldots,i_{\ell+1}\to j_{\ell+1}}| \ge |A_{i_1\to j_1,\ldots,i_{\ell}\to j_{\ell}}|(1-\frac{7}{n}), 
\end{equation}
and similarly for $B$. By a union bound,~\eqref{Eq:Boot1} and~\eqref{Eq:Boot2} imply that 
\begin{equation}\label{Eq:Boot3}
|A_{i_1 \rightarrow j_1,\ldots,i_{t}\to j_{t}}| \ge |A|(1-t\cdot \frac{7}{n}) \ge |A|(1-t\cdot \frac{7}{500t}) \ge 0.98|A|, 
\end{equation}
and similarly for $B$. In words, this means that both $A$ and $B$ are almost fully included in the $t$-umvirate $(S_n)_{i_1\to j_1,\ldots,i_t \to j_t}$. We now show that they are fully contained in that $t$-umvirate.

Assume that there exists $g \in B \setminus (S_n)_{i_1\to j_1,\ldots,i_t \to j_t}$. As $A$ and $B$ are cross $t$-intersecting, any $f \in A_{i_1\to j_1,\ldots,i_t \to j_t}$ agrees with $g$ on at least $t$ coordinates. By assumption, such $f,g$ can agree on at most $t-1$ among the coordinates $i_1,\ldots,i_t$, and hence, for each element $f \in A_{i_1\to j_1,\ldots,i_t \to j_t}$, there exists $i \not \in \{i_1,\ldots,i_t\}$ such that $f(i)=g(i)$. However, for any $g$, there are at least $\lfloor \frac{1}{e}(n-t)! \rfloor$ elements of $(S_n)_{i_1\to j_1,\ldots,i_t \to j_t}$ which do not agree with $g$ on any $i \not \in \{i_1,\ldots,i_t\}$ (as this number is clearly no less than the number of permutations on $n-t$ elements that have no fixed points), and hence,   
 \begin{equation*}
|A_{i_1\to j_1,\ldots,i_t \to j_t}| \le \lceil (1-\frac{1}{e})\cdot(n-t)! \rceil.
  \end{equation*}
Therefore,
\begin{align}\label{Eq:Boot4}
\begin{split}
|A|\cdot|B| &\le (\frac{100}{98})^2 \cdot |A_{i_1 \rightarrow j_1,\ldots,i_{t}\to j_{t}}||B_{i_1 \rightarrow j_1,\ldots,i_{t}\to j_{t}}| \\ &\le(\frac{100}{98})^2 \cdot (n-t)!\cdot \lceil (1-\frac{1}{e})\cdot(n-t)! \rceil < (n-t)!^2,    
\end{split}
\end{align}
in contradiction to the assumption $|A||B| \geq r(n,t)(n-t)!^2$. 

We assumed that there exists $g \in B \setminus (S_n)_{i_1\to j_1,\ldots,i_t \to j_t}$ and reached a contradiction. This implies that $B \subset (S_n)_{i_1\to j_1,\ldots,i_t \to j_t}$. By the same argument with the roles of $A,B$ reversed, we obtain $A \subset (S_n)_{i_1\to j_1,\ldots,i_t \to j_t}$. This completes the inductive proof of Assertion~(1).

\medskip

\paragraph{\emph{Proof of Assertion~(2).}} The first part of Assertion~(2) follows immediately from the proof of Assertion~(1), by considering $(t,n)$ such that $n \geq 2\lfloor c_0 t \rfloor$ and two cross $t$-intersecting families $A,B \subset S_n$ such that $|A||B| \geq \frac{3}{4} (n-t)!^2$ and applying to them Step~2 of the proof (namely, the `bootstrapping' step) verbatim. Note that as in the proof of Assertion~(1), we can apply Claim~\ref{Claim:Almost-in-dictatorship} $t$ times sequentially; the condition which corresponds to~\eqref{Eq:Boot1.5} is satisfied since $0.94 \cdot \frac{3}{4} > \frac{2}{3}$. 

The second part of Assertion~(2) (i.e., the case $A=B$) follows by noting that in this case,~\eqref{Eq:Boot4} can be replaced by 
\[
|A|\cdot|A| \le(\frac{100}{98})^2 \cdot (\lceil (1-\frac{1}{e})\cdot(n-t)! \rceil)^2 < (\frac{3}{4})^2 (n-t)!^2,
\]
which implies $|A| < \frac{3}{4} (n-t)!$, yielding a contradiction. This completes the proof of the proposition. 
\end{proof}

\bibliographystyle{plain}
\bibliography{refs}

\end{document}